\documentclass[oneside]{amsart}
\usepackage{microtype} 
\usepackage{my-settings}

\hypersetup{pdftitle={The combinatorics of associated Hermite polynomials}}

\graphicspath{{./}{figs/}}

\DeclareMathOperator{\lrm}{LRM}

\begin{document}
\begin{abstract}
  We develop a combinatorial model of the associated Hermite polynomials
  and their moments, and prove their orthogonality with a sign-reversing
  involution. We find combinatorial interpretations of the moments as
  complete matchings, connected complete matchings, oscillating
  tableaux, and rooted maps and show weight-preserving bijections
  between these objects. Several identities, linearization formulas, the
  moment generating function, and a second combinatorial model are also
  derived.
\end{abstract}

\title[Combinatorics of associated Hermites]
{The combinatorics of associated Hermite polynomials}
\author{Dan Drake}
\address{Department of Mathematical Sciences\\
         Korea Advanced Institute of Science and Technology\\
         Daejeon, South Korea}
\email{ddrake@member.ams.org}
\urladdr{http://math.kaist.ac.kr/~drake}
\date{\today}
\maketitle

The \emph{associated Hermite polynomials} are a sequence of orthogonal
polynomials considered by Askey and Wimp
in~\cite{askey.wimp:associated}, who analytically derived a number of
results about these polynomials. They are also treated in
\cite[Section~5.6]{ismail:classical}. In~\autoref{s:orthogonality} we
provide a combinatorial interpretation of these polynomials, their
moments, and describe an involution that proves the orthogonality and
$L^2$ norms of the polynomials with respect to those moments. Then
in~\autoref{s:linearizations} we shall describe several linearization
formulas involving associated Hermite polynomials. We finish with
weight-preserving bijections between a number of classes of
combinatorial objects whose generating functions all yield the moments
of the associated Hermites, and a second combinatorial model for the
polynomials.

We will assume that the reader is familiar with Viennot's general
combinatorial theory of orthogonal polynomials \cite{viennot:theorie,
viennot:combinatorial} and with the combinatorics of Hermite
polynomials; see \cite{azor.gillis.ea:combinatorial,
sainte-catherine.viennot:combinatorial,labelle.yeh:combinatorics} and
also \cite[\S II.6]{viennot:theorie}. In this paper we use $[n]$ to mean
the set of integers $1$ to $n$, inclusive, and write $[n] \sqcup [m]$
for the disjoint union of two such sets.

\section{Definition and orthogonality} \label{s:orthogonality}

The associated Hermite polynomials may be defined by shifting the
recurrence relation for the usual Hermite polynomials, which is
\begin{displaymath}
  H_{n+1}(x) = x H_n(x) - n H_{n-1}(x),
\end{displaymath}
to
\begin{equation}
  \label{e:associated-hermite-recurrence}
  H_{n+1}(x; c) = x H_n(x; c) - (n - 1 + c) H_{n-1}(x; c),
\end{equation}
with $H_0(x) = H_0(x; c) = 1$ and polynomials with negative indices
equal to zero. Askey and Wimp use a different normalization than we do;
one obtains our normalization from plugging $x/\sqrt{2}$ and $c-1$ into
their associated Hermites and dividing by $\left(\sqrt{2}\right)^n$.

The usual Hermite polynomial $H_{n+1}(x)$ is the generating function for
incomplete matchings of $[n+1]$, in which fixed points have weight $x$
and edges have weight $-1$; that combinatorial interpretation can be
derived from the recurrence relation as follows: the vertex $n+1$ may be
fixed with weight $x$, times the weight of all matchings on $[n]$; or we
may connect vertex $n+1$ to any of the $n$ vertices to its left, give
the edge weight $-1$ and multiply by all matchings on the $n-1$
remaining vertices.

For the associated Hermites, we'll build the matchings recursively as
described above and think of the parameter $c$ as meaning that one
special choice for the edge from $n+1$ will have weight $-c$. Two
natural choices are to make the special choice be the leftmost available
vertex, or the rightmost available vertex. Choosing the rightmost
available vertex happens to make the orthogonality involution easy to
prove, and yields the following result:

\begin{theorem}
  \label{t:assoc-hermite-comb-interp}
  The $n$th associated Hermite polynomial is the sum over weighted
  incomplete matchings $M$ of~$[n]$:
  \begin{equation}
    H_n(x; c) = \sum wt(M),
  \end{equation}
  in which fixed points have weight $x$, edges that nest no fixed points
  or edges and have no left crossings have weight $-c$, and all other
  edges have weight $-1$.
\end{theorem}

\begin{proof}
  We build the matching from right to left, and if at some point we add
  an edge and do not choose the rightmost available vertex, then that
  edge will nest a vertex, and when we come to that vertex, we will
  either leave it fixed (resulting in a fixed point underneath that
  edge), connect to another vertex underneath the edge (resulting in an
  edge nested by the original edge), or connect to a vertex to the left
  of the edge, resulting in a left crossing for the original edge. Any
  of these possibilities indicate that the rightmost vertex was not
  chosen, so edges for which none of those happen must have weight $-c$.
\end{proof}

\begin{figure}[ht]
  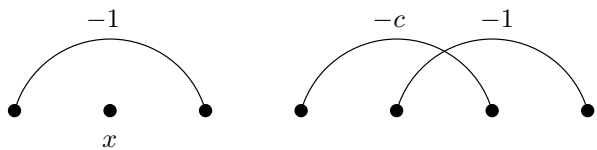
  \caption[A matching on $8$ points using the rightmost-choice
  weighting.]
  {A matching on $7$ points using the rightmost-choice weighting. This
  weighting is used throughout this paper.}
  \label{f:rightmost-choice-matching}
\end{figure}

An example of such a weighted matching is shown in
\autoref{f:rightmost-choice-matching}.

With nothing more than this model, we can easily explain an
``unexpected'' limit that Askey and Wimp derive
\cite[eq.~(5.9)]{askey.wimp:associated}. (In their paper, there is a
small typo: it should be $H_n(x \sqrt{2c}; c)$.) Using our
normalizations, the limit is
\begin{equation}
  \lim_{c \to \infty} c^{-n/2} H_n (x \sqrt{c}; c) = U_n(x),
  \label{e:assoc-hermite-to-cheby}
\end{equation}
where $\{U_n(x)\}_{n \ge 0}$ are the Chebyshev polynomials of the second
kind, also known as Fibonacci polynomials \cite[\S
II.1]{viennot:theorie}, \cite{sainte-catherine.viennot:combinatorial}.
$U_n(x)$ may be thought of as the generating function for incomplete
matchings on $n$ vertices in which edges always connect adjacent
vertices and have weight $-1$, and fixed points have weight $x$.

Using that combinatorial interpretation for $U_n(x)$ and the above
interpretation for $H_n(x;c)$, there is nothing unexpected about this
limit: take $H_n(x \sqrt{c}; c)$ and give each vertex, whether fixed or
incident to an edge, weight $1/\sqrt{c}$, so that $c^{-n/2}
H_n(x\sqrt{c}; c)$ is the generating function for incomplete matchings
with fixed points weighted $x$, and all edges weighted $-1/c$ except
those which nest no fixed points or edges, and have no left
crossings---such edges have weight $-1$. As $c$ goes to infinity, we
effectively restrict the generating function to matchings in which no
edge has weight $-1/c$; i.e., every edge nests no fixed points or
edges, and has no left crossings, so every edge must connect adjacent
vertices.

We want a linear functional $\L_c$ with respect to which the associated
Hermite polynomials are orthogonal. This linear functional is determined
by its moments $\L_c(x^n)$, which according to Viennot's general
combinatorial theory of orthogonal polynomials, can be expressed as a
sum over weighted Dyck paths in which a northeast edge has weight $1$
and a southeast edge leaving from height $j$ has weight $j-1+c$. There are
no Dyck paths of odd length, so the odd moments are zero. The first few
nonzero moments are
\begin{align*}
  \mu_0 &= 1,     &\mu_4 &= 2c^2 + c,\\
  \mu_2 &= c,     &\mu_6 &= 5c^3 + 7c^2 + 3c.
\end{align*}
Using the bijection from weighted Dyck paths to complete matchings from
\cite[\S II.6]{viennot:theorie}, we have two combinatorial
interpretations for the moments:

\begin{theorem}
  \label{t:moment-interp}
  The $n$th moment $\mu_n(c)$ of the associated Hermite polynomials is
  the generating function for complete matchings of $[n]$ weighted by
  either: (1) edges which are not nested by any other edge have weight
  $c$, and all other edges have weight $1$; or (2) edges with no right
  crossings have weight $c$ and all other edges have weight $1$.
\end{theorem}

\begin{figure}[ht]
  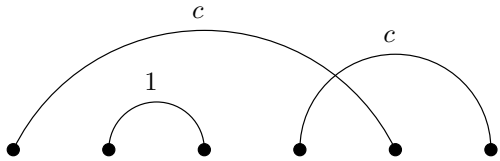
  \caption{A complete matching on $6$ points under the leftmost-choice
  weighting for the moments, in which nonnested edges have weight $c$
  and others have weight~$1$.}
  \label{f:leftmost-choice-moment}
\end{figure}

\begin{figure}[ht]
  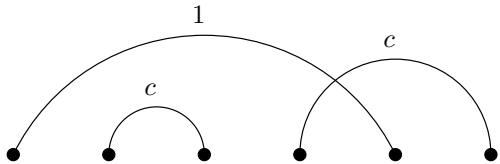
  \caption{The same complete matching under the rightmost-choice
  weighting, in which edges with no right crossing have weight $c$ and
  others have weight $1$. This is used in the orthogonality proof.}
  \label{f:rightmost-choice-moment}
\end{figure}

The two weightings correspond to giving weight $c$ to the leftmost and
rightmost choice, respectively, in the matchings. These interpretations
also explain why the odd moments are zero, since there are no complete
matchings on an odd number of vertices. For the proof of orthogonality,
we shall use the rightmost weighting; later we shall use the leftmost
weighting. Figures~\ref{f:leftmost-choice-moment} and
\ref{f:rightmost-choice-moment} show a matching using the two
weightings.


\subsection{Proof of orthogonality} \label{s:orthogonality-proof}

We wish to prove the following theorem in a combinatorial manner:

\begin{theorem}
  \label{t:assoc-hermite-orthogonality}
  The associated Hermite polynomials $H_n(x; c)$ are orthogonal with
  respect to the linear functional $\L_c$ with the above moments. They
  satisfy
  \begin{equation}
    \label{e:associated-hermite-integral}
    \L_c(H_n(x; c) H_m(x; c)) = \begin{cases}
      0               & n \ne m,\\
      \rising{c}{n} & n = m.
    \end{cases}
  \end{equation}
\end{theorem}

Here $\rising{a}{n}$ denotes the rising factorial $a (a+1)\cdots(a+n-1)$.

\begin{proof}
  The proof proceeds very similarly to the proof of orthogonality for
  usual Hermite polynomials. The product $H_n(x; c) H_m(x; c)$ is the
  generating function for pairs of matchings with, say, black edges,
  using the rightmost weighting. Applying $\L_c$ has the effect of
  putting a complete matching with the rightmost weighting with, say,
  green edges on the fixed points. We will use the phrase \emph{paired
  matching} to refer to such an object, with black homogeneous edges and
  arbitrary green edges, weighted as above. This is not standard
  terminology; it is only for our convenience.

  Using Theorems~\ref{t:assoc-hermite-comb-interp} and
  \ref{t:moment-interp}, the left side of
  \eqref{e:associated-hermite-integral} is the generating function
  for paired matchings, where black edges have weight $-c$ if they nest
  no edges, have no green crossings and no left black crossing;
  otherwise black edges have weight $-1$. Green edges have weight $c$ if
  they have no right green crossing and weight $1$ otherwise. See
  \autoref{f:paired-matching} for an example of such an object for $n =
  5$ and $m=3$.

\begin{figure}[ht]
  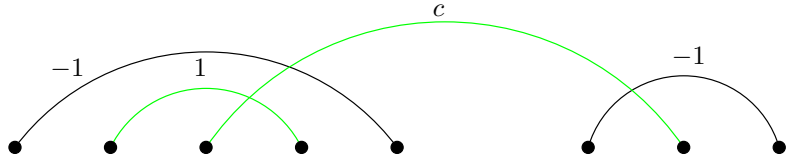
  \caption{A paired matching for $n=5$ and $m=3$.}
  \label{f:paired-matching}
\end{figure}

  We need an involution that shows the generating function for paired
  matchings equals zero when $n \ne m$, and equals $(c)_n$ otherwise.
  Assume that $n \ge m$ and put $[n]$ to the left of $[m]$. The involution
  is the very similar to that used in the combinatorial proof of
  orthogonality for usual Hermite polynomials: 
  \begin{quote}
    Find the leftmost homogeneous edge that nests no other edges and
    change its color.
  \end{quote}
  For example, in \autoref{f:paired-matching}, one would change the
  color of the leftmost green edge that connects vertices $2$ and $4$ to
  black. This operation is evidently an involution and will certainly
  change the sign; we need to verify that the weight of no other edge is
  affected by this change, and that if we change the color of an edge
  weighted $\pm c$ or $\pm 1$, the new edge has weight of $\mp c$ or
  $\mp 1$, respectively.
  
  We begin with the following observation: the leftmost homogeneous edge
  in $[n]$ that nests no edges can have no left crossing. We must check
  the four possibilities of color and weight for the edge whose color we
  flip:
  \begin{itemize}
    \item Edge is black, weight $-1$: the edge has no left crossing, and
    we've assumed the edge nests no edges, so if it has weight $-1$ it
    must have a right crossing by a green edge---so as green, it will
    have weight $+1$.
    \item Edge is black, weight $-c$: to get weight $-c$, the edge must
    in particular have no green crossing, and therefore as green, it
    will have weight $c$.
    \item Edge is green, weight $1$: the edge must have a right crossing
    by a green edge, so as black it will have weight $-1$.
    \item Edge is green, weight $c$: the edge nests no edges by
    assumption, and has no right green crossing. By our observation
    above, it has no left crossings, hence will be eligible for weight
    $-c$ as a black edge.
  \end{itemize}
  Thus the weight of the edge is preserved and the sign is reversed. We
  leave it to the reader to check that the weight and sign of no other
  edge is affected by this operation. It is only necessary to consider
  an edge that has a left crossing by the edge whose color changes.

  If $n > m$, there must be a homogeneous edge in $[n]$; in that case, the
  above involution has no fixed points, and we have proved that $H_n(x;c)$
  is orthogonal to $H_m(x;c)$.

  Now we shall prove that the $L^2$ norm of the associated Hermites is
  $\rising{c}{n}$ by interpreting the paired matchings as something
  whose generating function is known to be $\rising{c}{n}$: permutations
  weighted by left-to-right maxima. See \cite{medicis.viennot:moments,
  foata.strehl:combinatorics} for proofs of this fact in the context of
  Laguerre polynomials. (``Left-to-right maxima'' is ``\'el\'ements
  saillants inf\'erieurs gauches'' in French.) This bijection naturally
  generalizes the usual combinatorial proof that the $L^2$ norm of the
  Hermite polynomials is $n!$.

  First, apply the above involution to paired matchings with $n=m$; that
  involution will cancel all matchings with a homogeneous edge. To set
  up the bijection, begin with a matching on $[n] \sqcup [n]$ with no
  homogeneous edges. (Recall that $[n] \sqcup [n]$ means the disjoint
  union of $[n]$ with itself, or, what will work equally well, the
  ordinary union $\{1,\dots,n\} \union \{n+1,\dots,2n\}$.) Number the
  vertices as shown in \autoref{f:lrm-matching} and think of the right
  side as the domain, and the left side as the range.  A simple
  induction argument demonstrates that edges that get weight $c$
  correspond exactly to digits in the permutation that are left-to-right
  maxima.

  \begin{figure}[ht]
    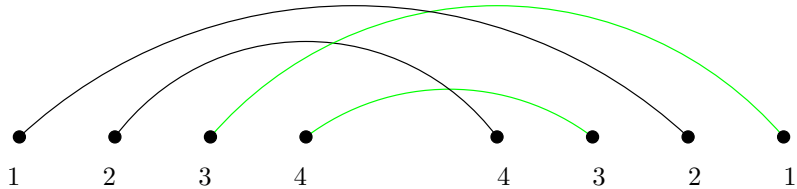
    \caption[The permutation $3142$ as a matching.]
    {The permutation $3142$ as a matching. The right side is the domain,
    the left is the range. The digits $3$ and $4$ are left-to-right
    maxima in the permutation, and indeed the green edges connecting to
    $3$ and $4$ on the left have weight $c$ under the associated Hermite
    moment weighting.}
    \label{f:lrm-matching}
  \end{figure}

  This bijection from the fixed points of the involution to
  permutations preserves weight, hence the $L^2$ norms of the
  associated Hermite polynomials are $\rising{c}{n}$. This completes
  the proof of \autoref{t:assoc-hermite-orthogonality}. 
\end{proof}

We also note that by \cite[Lemma~2.1]{foata.strehl:combinatorics}, the
$L^2$ norm can also be interpreted as the generating function for
permutations with cycles weighted by $c$.

\section{Linearizations} \label{s:linearizations}

In~\cite[theorem 3.1]{markett:linearization}, Markett shows that the
linearization coefficients in
\begin{equation}
  \label{e:linearization-of-assoc-hermites}
  H_N(x;c) H_M(x;c) = \sum_{j=0}^{\mathclap{\min(N,M)}} f(N,M,j)\ H_{N+M-2j}(x;c)
\end{equation}
are
\begin{equation}
  \label{e:assoc-linearization-coeffs}
  f(N,M,j) = \rising{N+M-2j+c}{j}\ \threeftwoatone{j-N}{j-M}{-j}{j-N-M-c+1}{1},
\end{equation}
where the ${}_{3}F_2$ notation indicates a hypergeometric function
evaluated at $x=1$. We can prove
\begin{theorem}
  \label{t:nonneg-poly-in-c}
  The linearization coefficients $f(N,M,j)$ of equation
  ~\eqref{e:assoc-linearization-coeffs} are polynomials in $c$ with
  nonnegative integer coefficients.
\end{theorem}

\begin{proof}
Take the rising factorial in front and reverse the order of
multiplication: it becomes $(-1)^j \rising{j-N-M-c+1}{j}$. We have two
$k!$ factors in the denominator of the ${}_3 F_2$; combine them with the
$\rising{j-N}{k}$ and $\rising{j-M}{k}$ in the numerator to get $(-1)^k
\binom{N-j}{k}$ and $(-1)^k \binom{M-j}{k}$. The $(-1)^k$ factors
cancel. Finally rewrite $\rising{-j}{k} = \rising{j-k+1}{k} (-1)^k$.

Put $(-1)^j \rising{j-N-M-c+1}{j}$ inside the sum. There is a factor of
$\rising{j-N-M-c+1}{k}$ in the denominator; those cancel and yield $(-1)^j
\rising{j-N-M-c+1+k}{j-k}$ in the numerator of the sum. Reverse the order
again and it turns into $(-1)^k \rising{N+M-2j+c}{j-k}$. This
$(-1)^k$ cancels with the earlier one from the $\rising{-j}{k}$.

The sum is now
\begin{displaymath}
  \sum_{k \ge 0}\binom{N-j}{k}\binom{M-j}{k} (j-k+1)_k
  \rising{N+M-2j+c}{j-k}.
\end{displaymath}
This is clearly a polynomial in $c$ with nonnegative coefficients. 
\end{proof}

Note that when $c=1$, the ${}_3 F_2$ of
\eqref{e:assoc-linearization-coeffs} sums by the Pfaff-Saalsch\"utz
identity to
\begin{displaymath}
  \frac{\rising{N+1-j}{j}\rising{M+1-j}{j}}{j!},
\end{displaymath}
and we recover the linearization coefficients for usual Hermite
polynomials; the expression above, after multiplying by $(N+M-2j)!$,
counts inhomogeneous matchings on $[N] \sqcup [M] \sqcup [N+M-2j]$, as
shown by de~ Sainte-Catherine and Viennot in
\cite{sainte-catherine.viennot:combinatorial} and, using different
methods, by Zeng in \cite{zeng:weighted}. A combinatorial interpretation
of the coefficients~ \eqref{e:assoc-linearization-coeffs}, refining the
results of \cite{sainte-catherine.viennot:combinatorial} and
\cite{zeng:weighted}, is quite desirable, but the problem is still open;
see \autoref{t:linearization-conj} below.

Since the linearization coefficients are known to be multiples of a
${}_3 F_2$ hypergeometric series, the best starting points for a
combinatorial interpretation seem to be \cite{andrews:identities,
nanjundiah:remark, andrews.bressoud:identities}; the first two papers
concern the usual Pfaff-Saalsch\"utz identity, the third features a
combinatorial proof of the $q$-Pfaff-Saalsch\"utz identity. It seems
very difficult to even prove, in analogy to the case for usual Hermite
polynomials, that $\L_c(H_N(x;c) H_M(x;c) H_{N+M-2j}(x;c))$ is the
generating function for inhomogeneous matchings on $[N] \sqcup [M]
\sqcup [N+M-2j]$.

In fact, using the ``nonnested'' weighting for the moments, the
generating function for inhomogeneous matchings doesn't even equal the
integral of three associated Hermites: $\L_c(H_2(x;c)^3) = c^3 + 4c^2 +
3c$, but the $8$ inhomogeneous matchings on $[2] \sqcup [2] \sqcup [2]$
have total weight $2c^3 + 4c^2 + 2c$. 

However, even if we use the ``no left crossing'' moments, it can be
shown that no involution that simply flips the color of an edge can work
with this model of the polynomials. For example, in $\L_c(H_1(x;c)
H_2(x;c) H_3(x;c))$, the matching $(1,2)(3,5)(4,6)$ has a homogeneous
edge---the one connecting $4$ and $6$---which as a black edge has weight
$-1$ and as a green edge has weight $+c$. That matching has only one
inhomogeneous edge, but changing its color does not preserve weight.

One might try different weightings for the polynomials and the moments.
We could reverse the matchings with the rightmost-choice weighting and
give edges with no \emph{left} crossing weight $c$. For the polynomials,
one could build them from left to right or right to left, and have
weight $-c$ given to the rightmost or leftmost choice. That yields two
moment weightings and four polynomial weightings, and counterexamples
like the one above are known for all eight combinations of weightings.

The order in which the sets of vertices are arranged is also important.
The integral $\L_c(H_3(x;c) \cdot H_3(x;c)\cdot  H_4(x;c))$ equals
$c(c+1)(c+2) (c+3)(c+8)$, but even if one considers only inhomogeneous
matchings, the three ways to arrange the sets of vertices yield three
different generating functions for inhomogeneous matchings with the
rightmost-choice moment weighting:
\begin{align*}
  [3] \sqcup [3] \sqcup [4]: & \ c(c+1)(c+2)(c+3)(c+8) \\
  [3] \sqcup [4] \sqcup [3]: & \ c(c+1)(c+2)(c^2 + 7c + 28)\\
  [4] \sqcup [3] \sqcup [3]: & \ c(c+1)(c+2)(c^2 + 8c + 27). 
\end{align*}
The nonnested weighting for the moments also fails in all three of these
cases: it gives $6 c(c+1)^2(c+2)^2$ for $[3] \sqcup [4] \sqcup [3]$ and
$3 c(c+1)(c+2)^2(c+3)$ for the other two. 
To get the correct answer, we had to order the sets of vertices in
weakly increasing order by size: $[3] \sqcup [3] \sqcup [4]$. 
This observation (and much computational evidence) leads us to conjecture
the following:

\begin{conjecture}
  \label{t:linearization-conj}
  Let $n_1,n_2,\dots,n_k$ be positive integers. The integral
  \[
  \L_c \left( \prod_{i=1}^k H_{n_i}(x;c) \right)
  \]
  is the generating function for inhomogeneous matchings on
  $\bigsqcup_{i=1}^k [n_i]$ in which the sets of vertices are arranged
  in weakly increasing order by size and the edges are weighted with the
  rightmost-choice moment weighting (so edges with no right crossing
  have weight $c$).
\end{conjecture}

\subsection{A mixed linearization formula}

In this section we will prove

\begin{theorem}
  \label{t:mixedlinearization}
  If $n \ge m-1$, then
  \begin{equation}
    \label{e:mixedlinearization}
    H_n(x;c) H_m(x) = \sum_k \binom{n-1+c}{k}\binom{m}{k}k!
    H_{n+m-2k}(x;c),
  \end{equation}
  where the sum runs from $0$ to $\min(m, \lfloor (n+m)/2 \rfloor)$.
\end{theorem}

\begin{proof}
  Fix $n$; we'll induct on $m$. For $m=0$ and $m=1$ the formula is a
  tautology and the recurrence relation, respectively. Assume that the
  formula works for some $m \le n$; multiply both sides of the formula
  by $x$ and use the recurrence:
  \begin{multline*}
    H_n(x;c) ( H_{m+1}(x) + m H_{m-1}(x) ) =\\
    \sum_k
    \binom{n-1+c}{k}\binom{m}{k}k! (H_{n+m+1-2k}(x;c) +
    (n+m-2k-1+c)H_{n+m-1-2k}(x;c)).
  \end{multline*}
  If we move the $m H_n(x;c) H_{m-1}(x)$ term over and use the
  induction hypothesis, we find that the coefficient of
  $H_{n+m+1-2k}(x;c)$ on the left side is
  \begin{displaymath}
    \binom{n-1+c}{k}\binom{m}{k}k! + (n+m-2k+1+c)
    \binom{n-1+c}{k-1}\binom{m}{k-1}(k-1)! -
    m\binom{n-1+c}{k-1}\binom{m-1}{k-1}(k-1)!
  \end{displaymath}
  which simplifies to
  \begin{displaymath}
    \binom{n-1+c}{k}\binom{m+1}{k} k!,
  \end{displaymath}
  exactly the coefficient we want. 
\end{proof}

One must be careful with that recurrence, though. If $k$ gets
too large the recurrence fails, because 
\begin{displaymath}
  xH_{-1}(x;c) = H_0(x;c) + (-2+c)H_{-2}(x;c)
\end{displaymath}
is \emph{false}. The induction argument works to go from $m=n$ to
$m=n+1$ because $xH_0 = H_1 + (-1+c) H_{-1}$, as long as one assumes
polynomials with negative indices are zero.

Is there is a natural extension or modification of the sum in
\eqref{e:mixedlinearization} when $n < m - 1$? The coefficient of
$H_{n+m-2k}(x;c)$ given in the sum is correct for $0 \le k \le n$
regardless of the relationship between $n$ and $m$ because of the
recurrence argument above, but there appears to be no particularly nice
or easy pattern to the coefficients of $H_{n+m-2k}(x;c)$ for $k > n$
when $n < m - 1$. 

\section{Associated Hermite moments and oscillating tableaux}
\label{s:moments-tableaux}

In this section we will describe a statistic on \emph{oscillating
tableaux}, also known as up-down tableaux, and a bijection between these
tableaux and complete matchings which is weight-preserving when using
the weight for associated Hermite moments. Oscillating tableaux were
described by Sundaram \cite{sundaram:cauchy}; see section $5$ of
\cite{chen.deng.ea:crossings} for discussion of their origins and the
bijection to complete matchings, and \cite{krattenthaler:growth} for an
extension of the results of \cite{chen.deng.ea:crossings} to fillings of
Ferrers diagrams.

\begin{figure}
  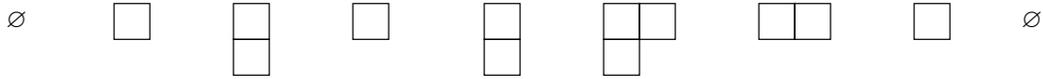
  \caption{An oscillating tableau that corresponds to the complete
  matching $(13)(26)(48)(57)$.}
  \label{f:tableau-to-matching} 
\end{figure}

Briefly, an oscillating tableau is a path in the Hasse diagram of the
Young lattice in which at each point one either moves up to a partition
that covers the current partition, or moves down to a partition covered
by the current partition. For our purposes, the path will always begin
and end with the empty shape. The \emph{length} of an oscillating
tableau is the number of edges in the path.
\autoref{f:tableau-to-matching} has an example of an oscillating tableau
of length~$8$.

In this section, we use \autoref{t:moment-interp}'s
``leftmost-available'' weighting of complete matchings, in which edges
that are not nested by other edges have weight $c$, and all other edges
have weight~$1$. 

\begin{figure}
  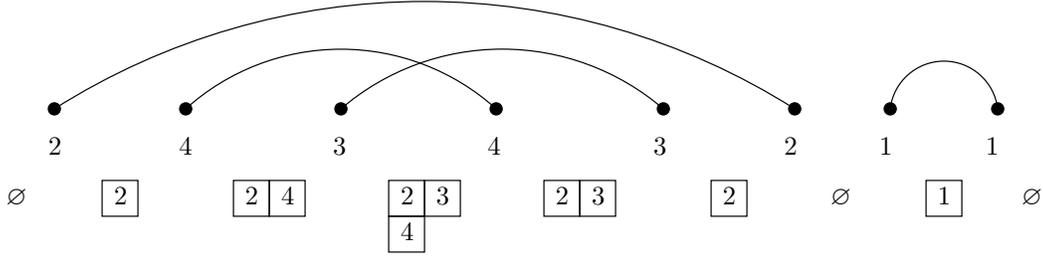
  \caption{A complete matching and the corresponding oscillating
  tableau. The numbers in the Ferrers shapes are not, strictly speaking,
  part of the oscillating tableau; they are only used in the bijection
  from the matching to the tableau.}
  \label{f:matching-to-tableau} 
\end{figure}

Roughly speaking, the bijection from complete matchings to oscillating
tableaux works by RSK-inserting numbers when edges start, and deleting
them when edges end. More precisely, given a complete matching, number
the edges from right to left as in \autoref{f:matching-to-tableau}.
(Equivalently, write the matching as a double occurrence word; see
\autoref{s:moments-maps}.) We will map this matching to a sequence of
Ferrers shapes. Begin with the empty Ferrers shape and read the matching
left to right. When edge~$j$ starts, RSK-insert a~$j$ into the tableau;
when edge~$j$ ends, delete the box containing~$j$. When done, erase the
numbers in the Ferrers shapes. \autoref{f:matching-to-tableau} has an
example.

There is a possible point of confusion here. A tableau in this context
is a path in the Hasse diagram of the Young lattice---a sequence of
Ferrers shapes. A standard Young tableau is a path that continually
moves up, and therefore it is simple to record the path with a single
Ferrers shape filled with numbers that strictly increase in rows and
columns. In \autoref{f:matching-to-tableau}, the Ferrers shapes are
written as Young tableaux, which is only for our convenience. The actual
image of the complete matching is the same sequence without the numbers
in the shapes. The reason for this is that RSK is a bijection, and one
can unbump numbers.

\autoref{f:tableau-to-matching} describes the inverse map from
tableaux to matchings. We read the sequence of Ferrers shapes from right
to left. Because of how we number the edges, the first box must have
a $1$ in it. In general, when the shape gets larger, we put the
next-largest number into the new box, because we've started a new edge.
The third shape from the right is $\begin{smallmatrix} 1 & 2 \\ 3
\end{smallmatrix}$, and the shape to its left must be
$\begin{smallmatrix} 1 \\ 3 \end{smallmatrix}$, because unbumping the
$2$ is the only way to produce the second shape. This oscillating
tableau corresponds to the matching $43412321$, using the
vertex-numbering scheme described above.

Let us weight oscillating tableaux with the following statistic: numbers
that appear in the first column only have weight $c$, and all other
numbers have weight $1$. That statistic is exactly what we need to prove
the following theorem.

\begin{theorem}
  \label{t:oscillating-tableaux-bijection}
  There is a weight-preserving bijection between oscillating tableaux of
  length $2n$ weighted with the above statistic and complete matchings
  weighted with the leftmost-available associated Hermite weighting.
\end{theorem}

We will use several preliminary results to prove this theorem.

\begin{lemma}
  \label{t:add-new-number-to-tableau}
  In an oscillating tableau, when a number is added to a shape, the
  corresponding edge is nested by all edges whose corresponding number
  in the shape is smaller, and has a left crossing from all edges whose
  corresponding number in the shape is bigger. Edges whose corresponding
  numbers never appear together in a shape neither nest nor cross one
  another.
\end{lemma}

For example, when we move from $\begin{smallmatrix} 2 & 4 \\ \phantom{4}
\end{smallmatrix}$ to $\begin{smallmatrix} 2 & 3 \\ 4 \end{smallmatrix}$
in \autoref{f:matching-to-tableau}, edge $3$ is nested by edge $2$
and has a left crossing from edge $4$. The proof of this is left to the
reader; it follows from the way the edges are numbered and in what order
we add numbers to the tableau.

The above lemma implies the following facts:

\begin{proposition}
  In an oscillating tableau, edges that get nested by other edges are
  exactly those whose number appears in the $2$nd, $3$rd, etc, column of
  a shape. Edges that have a right crossing are exactly those whose
  number appears in the $2$nd, $3$rd, etc row of a shape.
\end{proposition}

\begin{proof}[Proof of \autoref{t:oscillating-tableaux-bijection}]
The bijection between complete matchings and oscillating tableaux
clearly preserves weight: edges that do not get nested by another edge
must appear in the first column only. Note also that we could have used
the rightmost-available weighting from \autoref{t:moment-interp}; in
that case, we would have needed to make our statistic ``entries that
appear in the first row and stay there get weight $c$''.
\end{proof}

\section{Associated Hermite moments, rooted maps, and connected
matchings}
\label{s:moments-maps}

In addition to the weight-preserving bijection between associated
Hermite moments and oscillating tableaux of
\autoref{s:moments-tableaux}, there is a weight-preserving bijection
between associated Hermite moments and rooted maps. See
\cite{tutte:what, jackson.visentin:atlas} for introductions to maps,
which may be thought of as a graph along with an embedding into a
surface. A rooted map is a map in which one edge has been oriented.
There is an axiomatic construction of maps that makes it natural to
think of the edges in a map as pairs of half-edges or edge ends and we
will speak of edge ends in this section.

This connection is motivated by the normalizations used by
\cite{markett:linearization} and \cite{askey.wimp:associated}, both of
which use (rescaled versions of) $H_n(x;c+1)$. The first few moments for
those polynomials are
\begin{align*}
  \mu_0(c+1) &= 1,    &\mu_4(c+1) &= 2c^2 + 5c + 3,\\
  \mu_2(c+1) &= c+1,  &\mu_6(c+1) &= 5c^3 + 22c^2 + 32c + 15.
\end{align*}
On the one hand, those moments are simply the moments we've been working
with all along, except that now edges with no right crossing (or
nonnested edges, depending on which weighting one uses) may have weight
$1$ \emph{or} weight $c$. On the other hand, if those polynomials in $c$
are generating functions for some objects in which $c$, and not
$c+1$, is the weight, setting $c$ to $1$ gives us a count of how many
objects there are, which facilitates searching. Doing so yields 
\begin{displaymath}
  1, 1, 2, 10, 74, 706, 8162, 110410, 1708394,\dots
\end{displaymath}
which is
\href{http://www.research.att.com/~njas/sequences/A000698}{sequence
A000698} in \cite{sloane:on-line}. This sequence likely first appeared
in \cite{touchard:sur}; it counts connected matchings (see below).

In Table 1 of \cite[page 10]{arques.beraud:rooted}, Arqu\`es and
B\'eraud count rooted maps by number of edges and vertices; that table
also describes associated Hermite moments: the entry in the $n$th row
and $m$th column is the number of rooted maps with $n$ edges and $m$
vertices, and is also the coefficient of $c^{m-1}$ in $\mu_{2n}(c+1)$. We
will weight each vertex in such a map by $c$ \emph{except} the vertex at
the head of the root edge, and use the bijection between rooted maps in
orientable surfaces and connected matchings found in the work of Ossona
de~Mendez and Rosenstiehl \cite{ossona-de-mendez.rosenstiehl:encoding,
ossona-de-mendez.rosenstiehl:connected}. A \emph{connected matching} on $2n$
vertices is one in which all vertices except $1$ and $2n$ are nested by
an edge. Equivalently, one can write a matching as a \emph{double
occurrence word} in the letters $1, 2,\dots, n$ where each letter
appears exactly twice; then a matching is connected if the corresponding
double occurrence word cannot be written as the concatenation of two
double occurrence words.

A double occurrence word corresponds to the vertex-numbering scheme used
in \autoref{s:moments-tableaux}. We shall weight connected matchings by
giving weight $c$ to all nonnested edges \emph{except} the edge
containing vertex $1$. Then we have

\begin{theorem}
  \label{t:maps-to-connected-matchings}
  The function given in \cite{ossona-de-mendez.rosenstiehl:encoding} and
  \cite{ossona-de-mendez.rosenstiehl:connected} is a weight-preserving
  bijection from rooted maps in orientable surfaces with $k$ vertices
  and $n$ edges to connected matchings on $2n+2$ vertices of weight
  $c^{k-1}$.
\end{theorem}

\begin{proof}
  The idea of the bijection is this: number the edges in the rooted map,
  add a new loop at the vertex adjacent to the root, then build a double
  occurrence word by visiting each vertex and adding the edge numbers
  adjacent to the vertex to the word.

  The bijection is weight-preserving because when deciding the next
  vertex to visit, the algorithm chooses the vertex in the rooted map
  corresponding to the leftmost unattached vertex in the
  partially-constructed matching. As we add edge ends to the list, we
  will add a new edge to the matching that contains that leftmost
  unattached vertex. No edge can then nest the newly created edge, so
  every visit to a new vertex in the rooted map results in exactly one
  nonnested edge in the matching.
\end{proof}

\begin{figure}[h]
  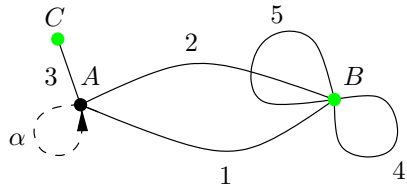
  \caption{A rooted map to which we'll apply the bijection to connected
  matchings. Green vertices have weight $c$. We have already added the
  extra edge, labeled $\alpha$; the original root was the end of edge
  $1$ incident with vertex $A$.}
  \label{f:rooted-map}
\end{figure}

\autoref{f:rooted-map} shows an example of the bijection. We will color
green the vertices of weight $c$ in the rooted map and the edges of
weight $c$ in the connected matching. We start at the head of edge
$\alpha$ and read counterclockwise around vertex $A$; our double
occurrence word begins with
\begin{displaymath}
  \alpha \ 1 \ 2 \ 3\ \alpha.
\end{displaymath}
We have visited both ends of $\alpha$, so we move to the unvisited end
of edge $1$, go around vertex $B$ and add $4\ 4\ 5\ 2\ 5\ 1$ to the
word, which is now
\begin{displaymath}
  \alpha \ 1 \ 2 \ 3\ \alpha\ 4\ 4\ 5\ 2\ 5\ 1.
\end{displaymath}
Now move to the unvisited end of edge $3$ and do the same thing; we just
append $3$ to the word. We end up with
\begin{displaymath}
  \alpha \ 1 \ 2 \ 3\ \alpha\ 4\ 4\ 5\ 2\ 5\ 1\ 3,
\end{displaymath}
which is double-occurrence word for the connected matching $(1, 5)(2,
11)(3, 9)(4, 12)(6, 7)(8, 10)$ where the edges $(2,11)$ and $(4,12)$
have weight $c$ because edges $1$ and $3$ in the rooted map were the
edges along which we first visited vertices $B$ and $C$, and $1$ and $3$
appeared in the double-occurrence word n positions $2$ and $11$, and $4$
and $12$ respectively.

Now we need another weight-preserving bijection, this time from weighted
connected matchings to one of our original definitions for $\mu_n(c+1)$,
the moments of associated Hermite polynomials. We will demonstrate such
a bijection to the moments weighted with the leftmost-available
weighting of \autoref{t:moment-interp}, in which nonnested edges are may
have weight $1$ or $c$. Call the edge containing vertex~$1$ the ``fake
edge''.

\begin{figure}[h]
  \includegraphics[scale=.75]{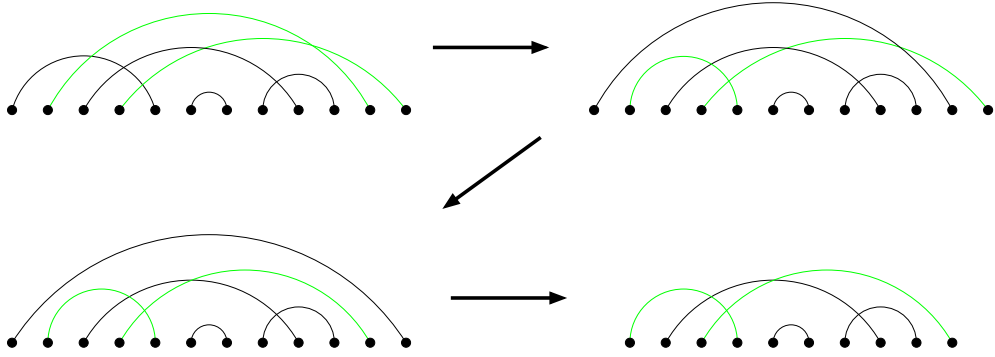}
  \caption{The steps of the tail-swapping bijection applied to the
  connected matching corresponding to the rooted map in
  \autoref{f:rooted-map}; the result is a complete matching (in the
  lower right) in which nonnested edges are eligible for weight $c$.}
  \label{f:conn-matching-tailswap}
\end{figure}

The bijection works as follows: If the fake edge has no crossings,
remove it; the remaining matching on $2n$ vertices, of weight $1$, is
the result of the bijection. Otherwise, swap the tails of the fake edge
and that edge crossing the fake edge which has the leftmost endpoint. That
crossing edge must have weight $c$; give the new edge, which is now
nested by the fake edge, weight $c$ also. Continue this tail-swapping
process with the fake edge until the fake edge has no crossings, then
remove it. An example is shown in \autoref{f:conn-matching-tailswap}.

This map is a bijection because it can be reversed: given such a
weighted matching on $2n$ vertices, add a new edge that nests the entire
matching, and swap tails with the green edges (those of weight $c$) from
right to left. Observe that the green edges in the connected
matching---which are nonnested---will end up nonnested after the
tail-swapping bijection, and vice versa, so this bijection is
weight-preserving. Note that in the example of
\autoref{f:conn-matching-tailswap} and \autoref{f:big-example-table},
the connected matching corresponded to a complete matching which was
also connected. Of course this does not always happen: the connected
matching $(1,5)(2,4)(3,8)(6,7)$ corresponds to the unconnected complete
matching $(1,3)(2,4)(5,6)$ under this bijection.

\autoref{t:maps-to-connected-matchings} established that the
generating functions for rooted maps and connected matchings are the
same; that theorem, together with the bijection between connected
matchings and arbitrary complete matchings, provides a proof of the
following theorem.

\begin{theorem}
  \label{t:maps-connected-complete-matchings}
  The generating functions for rooted maps with $n$ edges, connected
  matchings on $2n+2$ vertices, and complete matchings on $2n$ vertices
  all equal the moment $\mu_{2n}(c+1)$ of the associated Hermite polynomials.
\end{theorem}

\subsection{The moment generating function} Let $f(t; c)$ be the ordinary
generating function for the moments of the associated Hermite
polynomials:
\begin{displaymath}
  f(t; c) := \sum \mu_n(c) t^n.
\end{displaymath}
With the results of this section, we see that a continued fraction for
$f(t;c)$ is implicit in \cite[Theorem~3]{arques.beraud:rooted}: their function
$M(y,z)$ counts rooted maps with the exponent of $y$ counting the number
of vertices, and the exponent of $z$ counting the number of edges. We
know that $\mu_{2n}(c+1)$ is the generating function for rooted maps with
$n$ edges, in which all vertices except one get weight $c$, which means 
\begin{equation}
  f(t; c+1) = \frac{M(c, t^2)}{c} = 
  \cfrac{1}{1 - \cfrac{(c+1)t^2}{1 - \cfrac{(c+2)t^2}{1 -
  \cfrac{(c+3)t^2}{1 -\cdots }}}}.
  \label{e:moment-gf}
\end{equation}
This continued fraction can also be obtained with the method of \cite[p.
V-4]{viennot:theorie}, where Viennot shows a continued fraction
expansion for the moment generating function for any set of orthogonal
polynomials where the recurrence coefficients are known.

In the last two sections, we've shown bijections between the moments of
the associated Hermites, connected matchings, rooted maps and
oscillating tableaux. We summarize these correspondences by going all
the way from a rooted map, to a connected matching, to a regular
complete matching, to an oscillating tableau in
\autoref{f:big-example-table}.

\begin{table}
  \caption{A rooted map, a connected matching, a complete matching, and
  an oscillating tableau, all of weight $c^5$, that correspond to each
  other using the weight-preserving bijections of this paper. In the
  tableau, we have only colored the first box that corresponds to a
  number which gets weight $c$.}
  \label{f:big-example-table}
  \begin{tabular}{|c|p{1.4in}|}
    \hline
    Object & \makebox[1.4in]{What gets weight $c$} \\
    \hline
    \rule{0in}{.9in}
    \includegraphics[scale=.5]{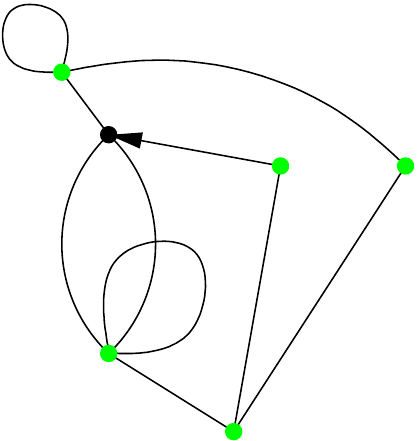} &
    \parbox[b][.9in][c]{1.4in}{Vertices not adjacent to head of root edge.} \\ 
    \hline
    \rule{0in}{1.25in}
    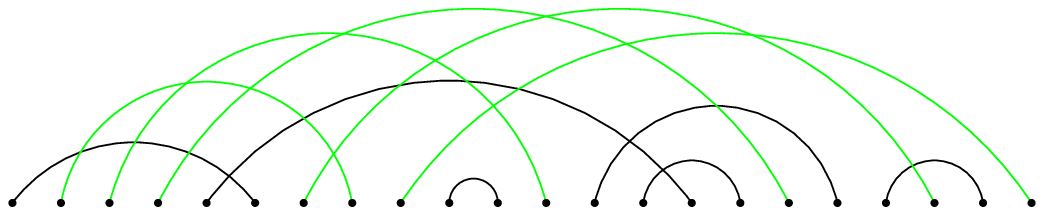 &
    \parbox[b][1.25in][c]{1.4in}{Non-nested edges except edge containing
    vertex $1$.} \\
    \hline
    \rule{0in}{1.20in}
    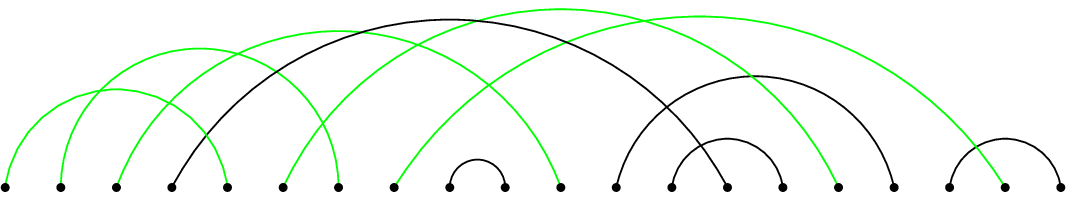 &
    \parbox[b][1.20in][c]{1.4in}{Non-nested edges may have weight $1$ or $c$.}\\
    \hline
    \rule{0in}{1.15in}
    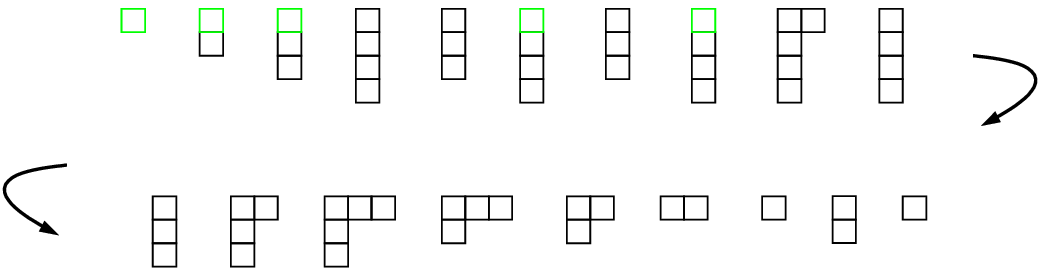 &
    \parbox[b][1.15in][c]{1.4in}{Numbers that appear in first column may
    have weight $1$ or $c$.} \\
    \hline
  \end{tabular}
\end{table}

\subsection{A second model for associated Hermite polynomials}
\label{s:fake-edge-model-polys}

The above discussion of connected matchings meshes nicely with a second
combinatorial model of the associated Hermites, which is motivated by
identity~ \eqref{e:identity1} below. The key features of this second
model are very similar to those of the connected matching model for the
moments: we are using $c+1$ but there are no choices for the weights of
parts of the matching, and the resulting matching is connected. The
identity is found in Askey and Wimp~\cite[equation~
(4.18)]{askey.wimp:associated} and we present a combinatorial proof.

\begin{theorem}
  \label{t:identity1}
  The associated Hermites may be written as a sum of usual Hermite
  polynomials:
  \begin{equation}
    H_n(x;c+1) = \sum_{k \ge 0} (-1)^k \rising{c}{k} \binom{n-k}{k} H_{n-2k}(x).
    \label{e:identity1}
  \end{equation}
\end{theorem}

We will need two lemmas to prove \autoref{t:identity1}.

\begin{lemma}
  \label{t:claim1}
  $(-1)^k \rising{c}{k}$ is the generating function for complete
  matchings on $2k$ vertices, with the $c+1$ associated Hermite polynomial
  weighting, such that all edges of weight $-1$ have a left crossing by
  an edge of weight $-c$. Furthermore, in such matchings there are
  exactly $k$ ``slots'' available underneath the edges weighted $-c$
  where one could place the left endpoint of a new edge of weight $-1$,
  and only one ``slot'' available for the left endpoint of a new edge of
  weight $-c$.
\end{lemma}

\begin{figure}
  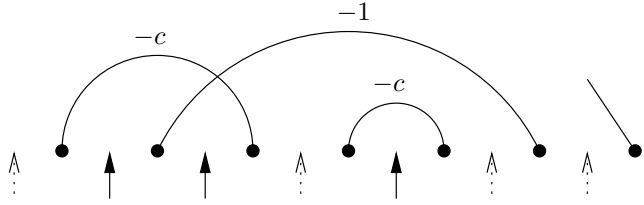
  \caption{A matching on $6$ vertices of the type described by
  \autoref{t:claim1}. If the new edge on the right is to have weight
  $-1$ and satisfy the conditions of the lemma, it must connect to a new
  vertex in one of the three available slots, indicated by the solid
  arrows.}
  \label{f:-1crossedby-c}
\end{figure}

\autoref{f:-1crossedby-c} shows an example of such a configuration.

\begin{proof}
  The proof goes by induction. The base cases are clear, and if true for
  some $k$, given any configuration for that $k$, we can either:
  \begin{itemize}
    \item add a new edge connecting vertices $2k+1$ and $2k+2$ which has
      weight $-c$, and hence we multiply the generating function for
      $2k$ vertices by $-c$ and add a new slot, or
    \item add a new edge from the rightmost vertex and put its
      left endpoint in any one of the $k$ ``slots'' underneath one of
      the $-c$ edges. Such an edge must have weight $-1$, and there are
      $k$ ways to place this edge, hence we effectively multiply the
      generating function by $k$, and since we put a new edge into one
      of the $k$ slots, there are now $k+1$ slots available below edges
      weighted $-c$.
  \end{itemize}
  See \autoref{f:-1crossedby-c} for an example of case 2. Altogether
  we've multiplied $(-1)^k \rising{c}{k}$, the generating function for
  $2k$ vertices, by $-(c+k)$, so the lemma is true by induction.
\end{proof}

\begin{lemma}
  \label{t:claim2}
  For such a configuration on $2k$ vertices as described in
  \autoref{t:claim1}, there are $k+1$ places where the left endpoint of
  one or more green edges of weight $1$ could be placed without
  affecting the weight of the configuration.
\end{lemma}

\begin{proof}
  Induction again. The green edges cannot cross the $-c$ edges. For
  example, in \autoref{f:-1crossedby-c}, there are four places where one
  could place such an edge, indicated by the dotted arrows.
\end{proof}

\begin{proof}[Proof of \autoref{t:identity1}]
  Since $H_{n}(x;c)$ is an even or odd polynomial if $n$ is even or odd,
  respectively, we can certainly write
  \begin{equation}
    H_n(x;c+1) = \sum_{k \ge 0} a_{nk} H_{n-2k}(x)
    \label{e:identity1pf}
  \end{equation}
  for some coefficients $a_{nk}$. We show that those coefficients equal
  $(-1)^k \rising{c}{k} \binom{n-k}{k}$. Fix $k$ between $0$ and $n/2$,
  multiply both sides by $H_{n-2k}(x)$, and apply the usual Hermite
  linear functional $\L_1$. On the right side, we use orthogonality and
  equation~\eqref{e:identity1pf} becomes
  \begin{displaymath}
    \L_1(H_n(x;c+1) H_{n-2k}(x)) = a_{nk} (n-2k)!.
  \end{displaymath}
  Thinking of the left side as paired matchings on $[n]$ and $[n-2k]$
  with black edges of weight $-1$ and $-c$ as appropriate, and green
  edges all of weight $1$, we may apply the following involution: find
  the leftmost homogeneous edge of weight $\pm 1$ in $[n]$ or $[n-2k]$
  and flip its color, \emph{unless} that edge has a left crossing with
  an edge of weight $-c$. Swapping the colors on such edges does not
  preserve the weight of the paired matching.

  \autoref{t:claim1} tells us the generating function of the
  configurations of edges that remain in $[n]$ after applying the
  involution; \autoref{t:claim2} tells us that such configurations may
  be viewed as consisting of $k$ ``chunks'' of vertices. Placing the
  green edges into those chunks is equivalent to forming a weak
  composition of $k$ into $n-2k+1$ parts; there are $\binom{n-k}{k}$
  such compositions, and having chosen where the $n-2k$ edges in $[n]$
  start, we can choose their endpoints in $[n-2k]$ in $(n-2k)!$ ways.
  Together we have
  \begin{displaymath}
    (-1)^k \rising{c}{k} \binom{n-k}{k}(n-2k)! =
    \L_1(H_n(x;c) H_{n-2k}(x)) = a_{nk} (n-2k)!
  \end{displaymath}
  which proves the identity of \autoref{t:identity1}.
\end{proof}

The above proof relies crucially on being able to give weight $-1$ or
$-c$ to edges; if we used $H_n(x;c)$, the above involution would not
cancel as many edges, and we would need to replace \autoref{t:claim1}
with something more complicated in order to handle the $\rising{c-1}{k}$
factor.

Our first model for the associated Hermite polynomials
(\autoref{t:assoc-hermite-comb-interp}) involved incomplete matchings on
$n$ vertices; the above identity motivates the following model for
$H_{n}(x;c+1)$ as matchings on $n+2$ vertices.

\begin{theorem}
  \label{t:assoc-hermite-interp-fake-edge}
  The associated Hermite polynomial $H_n(x;c+1)$ is the generating
  function for certain connected incomplete matchings on $n+2$
  vertices with the following weights:
  \begin{itemize}
    \item The edge containing vertex $1$ has weight $1$. Call
      this edge the ``fake edge''.
    \item Fixed points have weight $x$.
    \item Non-nested edges (except the fake edge) have
      weight $-c$.
    \item Nested edges have weight $-1$.
  \end{itemize}
  In such matchings, fixed points must be nested by the fake edge. All
  edges other than the fake edge must either cross or be nested by the
  fake edge.
\end{theorem}

\begin{figure}
  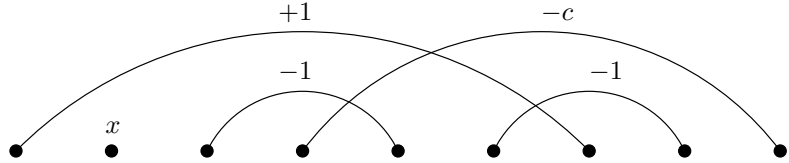
  \caption{A matching on $9$ vertices that contributes $cx$ to
  $H_7(x;c+1)$ using the combinatorial interpretation of
  \autoref{t:assoc-hermite-interp-fake-edge}. Note the ``fake edge'' of
  weight~$+1$.}
  \label{f:fake-edge-matching}
\end{figure}

An example of such a matching for $H_7(x;c+1)$ is shown in
\autoref{f:fake-edge-matching}. It is clear that the requirement for
nesting and crossing the fake edge yields a connected matching. Note
that the connected matching moments of \autoref{s:moments-maps} also
have a fake edge.

\begin{proof}[First proof]
  Consider the $k$th term in the sum~\eqref{e:identity1}:
  \begin{displaymath}
    (-1)^k \rising{c}{k} \binom{n-k}{k} H_{n-2k}(x).
  \end{displaymath}
  Begin with the fake edge and put $k$ vertices to the right of it. Put
  the remaining $n-k$ vertices underneath the fake edge and choose $k$
  of them to connect with the edges that will come from the $k$ vertices
  on the right of the fake edge; that accounts for the binomial
  coefficient. On the remaining $n-2k$ vertices underneath the fake
  edge, we put a regular Hermite-style matching; all the edges will have
  weight $-1$ since they are nested by the fake edge.

  The last thing to do is account for the $k$ edges that come from the
  right of the fake edge and show that they contribute $(-1)^k
  \rising{c}{k}$. According to \autoref{t:nonnested-edges-inc-subseq},
  the generating function for such a configuration with edges of weight
  $+1$ and $+c$ is $\rising{c+1}{k-1}$, but in our subset, the leftmost
  edge also gets weight $c$, so the correct factor is
  $\rising{c+1}{k-1}\cdot c = \rising{c}{k}$. Also, we must correct for
  the signs: our edges have weight $-1$ and $-c$, so we multiply by
  $(-1)^k$.
\end{proof}

\begin{proof}[Second proof]
  Verify that the generating function described in the theorem satisfies
  the three-term recurrence for the associated Hermites
  ~\eqref{e:associated-hermite-recurrence}. We proceed very much like
  the usual combinatorial proof of the recurrence relation for Hermite
  polynomials: any such restricted matching on $n+3$ vertices may be
  obtained by placing the fake edge and considering the rightmost vertex
  nested by the fake edge. There are three possibilities: one, we can
  leave that vertex fixed, and fill in the remaining $n+2$ vertices with
  any restricted matching; two, we can add an edge from that vertex to
  the very rightmost vertex, and fill in the remaining $n+1$ vertices
  with any restricted matching; three, we can attach that vertex to any
  vertex \emph{except} the rightmost vertex and fill in the remaining
  $n+1$ vertices as before. The first case contributes $x$ times the
  generating function for $n+2$ vertices. The second cases contributes
  $-c$ times the generating function for $n+1$ vertices, since that new
  edge cannot be nested, and it will not nest any of the other edges. In
  the third case, there are $n$ vertices to choose from and all of them
  will result in a nested edge of weight $-1$, so we add $-n$ times the
  generating function for $n+1$ vertices. This exposition is simply
  another way of stating~\eqref{e:associated-hermite-recurrence}:
  \begin{displaymath}
    H_{n+1}(x;c+1) = x H_n(x;c+1) - (n+c) H_{n-1}(x;c+1). \qedhere
  \end{displaymath}
\end{proof}

The following lemma was used in the first proof of
\autoref{t:assoc-hermite-interp-fake-edge}. It may be proved by
induction, similar to \autoref{t:claim1} and
\autoref{t:assoc-hermite-orthogonality}.

\begin{lemma}
  \label{t:nonnested-edges-inc-subseq}
  The generating function for complete matchings on $2n$ vertices in
  which all edges go from the ``left $n$'' vertices to the ``right $n$''
  vertices , with all nonnested edges having weight $c$ except the edge
  containing the leftmost vertex, is $\rising{c+1}{n-1}$.

  There is a weight-preserving bijection between such matchings and
  permutations $\pi$ of $[n]$ weighted by $c^{\,\lrm(\pi)-1}$ where
  $\lrm(\pi)$ is the number of left-to-right-maxima of the
  permutation.
\end{lemma}

At this point, we have a combinatorial interpretation for both the
associated Hermite polynomials
(\autoref{t:assoc-hermite-interp-fake-edge}) and their moments
(\autoref{t:maps-connected-complete-matchings}) in terms of connected
matchings with a fake edge; the natural thing to do is combine these to
get another proof of orthogonality. This will be quite difficult because
it is not at all obvious how to combine a pair of matchings for the
polynomials and a matching for the moments to get a paired matching; one
would have two fake edges from the polynomials and would need to somehow
incorporate the fake edge from the moments into that configuration.
However, it is interesting to note that the above theorem tells us how
we would derive the $L^2$ norm using such a setup: $H_n(x;c)^2$ would be
a pair of matchings on $2n+4$ vertices, but because of the extra fake
edge mentioned above, after canceling all homogeneous edges we would
effectively get complete matchings on $2n+2$ vertices in which all the
edges go from the left $n+1$ vertices to the right $n+1$. The generating
function for such a configuration, according to
\autoref{t:nonnested-edges-inc-subseq}, is $\rising{c+1}{n}$, which
agrees with the known $L^2$ norm for the associated Hermites at $c+1$.

\section{Unanswered questions and future directions}

We have taken the basic combinatorial model in \autoref{s:orthogonality}
for associated Hermite polynomials and their moments and gone in two
directions: to oscillating tableaux, and to rooted maps. The appeal of
oscillating tableaux is in the recent flurry of work on $k$-crossings
and $k$-nestings in matchings and set partitions; see
\cite{chen.deng.ea:crossings, krattenthaler:growth, mier:k-noncrossing,
klazar:identities, kasraoui.zeng:distribution, jelnek:dyck}. The moments
of Charlier polynomials are generating functions for set partitions and
it seems likely that some of this work could be used to treat the
associated Charlier polynomials.

Observe that in the connected matchings, the rooted maps, and in the
second combinatorial model for the associate Hermite polynomials of
\autoref{t:assoc-hermite-interp-fake-edge}, each model has some sort of
``fake edge''. Combining the models for the moments and polynomials
which both involve connected matchings would be interesting, but this
has not yet shown promise. A major problem is that each incomplete
matching for the polynomial is weighted by $x$ to the number of fixed
points---say there are $2k$ fixed points---but the corresponding
matchings are matchings on $2k+2$ vertices. It is not clear how to
combine these two objects in a geometric or graph-theoretical way that
allows a natural and easy proof of orthogonality.

Using rooted maps holds promise, though: Ossona de~Mendez and
Rosenstiehl have generalized the bijection between connected matchings
and rooted maps to a bijection between permutations and hypermaps
\cite{ossona-de-mendez.rosenstiehl:transitivity,
ossona-de-mendez.rosenstiehl:connected}. This suggests an intriguing
connection to Laguerre polynomials since hypermaps are built out of
permutations in the same way that maps are built out of complete
matchings. The paper of Askey and Wimp \cite{askey.wimp:associated}
which inspired this work devotes much more attention to the associated
Laguerres than to Hermites---about two thirds of the article. It is
natural, then, to work out a corresponding combinatorial treatment of
those polynomials, especially given the connections between rooted maps
and hypermaps. There is also the work of Ismail et al.\ 
\cite{ismail.letessier.ea:linear} who work with the associated Laguerres
as birth and death processes---there has been work on birth and death
processes and lattice paths \cite{flajolet.guillemin:formal} which
suggests another avenue for a combinatorial theory of those polynomials.

\section{Acknowledgements}

This work is based on part of the author's doctoral thesis, completed at
the University of Minnesota under the direction of Dennis Stanton. The
author thanks Professor Stanton for his assistance and patience and the
University of Minnesota math department for their support. This work was
presented at \href{http://fpsac.cn/}{FPSAC 2007} and the author thanks
the FPSAC referees for their careful reading and helpful comments.
Thanks also to Bill Chen and Jang Soo Kim, who helped correct some minor
errors.

\bibliographystyle{amsalphaurl} \bibliography{bib}
\end{document}